\documentclass[11pt]{article}

\usepackage{amssymb,amsthm, graphicx, stmaryrd, fancyhdr}
\usepackage{amsfonts}
\usepackage{amsmath}
\usepackage[english,german]{babel}
\usepackage{epsfig,epstopdf}
\usepackage{subfigure}
\usepackage{latexsym}
\usepackage{color}

\usepackage{amsmath,amssymb,epsfig, epstopdf}

\newtheorem{theorem}{Theorem}[section]
\newtheorem{definition}[theorem]{Definition}

\newtheorem{lemma}[theorem]{Lemma}

\newtheorem{remark}[theorem]{Remark}

\newtheorem{assumption}[theorem]{Assumption}

\def\R{\mathbb{R}}

\def\P{\mathbb{P}}

\def\R{\mathbb{R}}

\def\P{\mathbb{P}}

 \normalsize

\begin{document}
\selectlanguage{english}

 \title{On Investment-Consumption  with Regime-Switching \thanks{  Work supported by NSERC grants 371653-09, 88051 and MITACS grants 5-26761, 30354 and the Natural Science Foundation of China (10901086).}}

\author{\normalsize    Traian A.~Pirvu \\[8pt]
        \small Dept of Mathematics \& Statistics\\
        \small McMaster University \\
        \small 1280 Main Street West \\
        \small Hamilton, ON, L8S 4K1\\
        \small tpirvu@math.mcmaster.ca
        \and
        \normalsize Huayue Zhang \\[8pt]
        \small Dept of Finance\\
        \small  Nankai University\\
        \small 94 Weijin Road \\
        \small Tianjin, China, 300071 \\
        \small  hyzhang69@nankai.edu.cn
\vspace*{0.8cm}}

\maketitle

\noindent {\bf Abstract.}
 In a continuous time stochastic economy, this paper considers the problem of consumption and investment in a financial market in which the representative investor exhibits a change in the discount rate. The investment opportunities are a stock and a riskless account. The market coefficients and discount factor switches according to a finite state Markov chain. The change in the discount rate leads to time inconsistencies of the investor's decisions. The randomness in our model is driven by a Brownian motion and Markov
 chain. Following \cite{EkePir} we introduce and characterize the equilibrium policies for  power utility functions. Moreover, they are computed in closed form for logarithmic utility function. We show that a higher discount rate leads to a higher equilibrium consumption rate. Numerical experiments show the effect of both time preference and risk aversion on
 the equilibrium policies.
\vspace{1cm}


\begin{quote}

\end{quote}

\begin{flushleft}

\end{flushleft}

\setcounter{equation}{0}
\section{Introduction}

Dynamic asset allocation in a stochastic paradigm received
a lot of scrutiny lately. The first papers in this area are
\cite{Mer69} and \cite{Mer71}. Many works then followed, most
of them assuming an exponential discount rate. \cite{EkePir} has given an overview of the literature in
the context of Merton portfolio management problem with exponential
discounting.

The issue of discounting was the subject of many studies in financial economics.
Several papers stepped away from the exponential discounting modeling, and based
on empirical and experimental evidence proposed different discount models. They can
be organised in two classes: exogenous discount rates and  endogenous  discount rates.
In the first class the most well known example is the hyperbolic discounting. This type
discounts near future more heavily than distant future which is in accordance with the
experimental findings.  \cite{EkePir} and  \cite{EkePirMbo} discuss about this class of discounting.

The concept of endogenous time preference was developed by \cite{Koop} in a discrete
time formulation. \cite{Uza} considered the continuous time version which was later extended
by \cite{Eps}. The class of discount rates emerged in response to the following two observed
phenomena: ``decreasing marginal impatience'' DMI and  ``increasing marginal impatience''
IMI. DMI means that the lower the level of consumption, the more heavily an agent discount the future.
IMI is just the opposite: the higher the level of consumption, the more heavily an agent discount the future.
Some papers support DMI, e.g. \cite{Das}, others advocate for IMI, \cite{Tak}.

In this paper we consider a regime switching model for the financial market. This modeling is
consistent with some cyclicality observed in financial markets. Many papers considered these
types of markets for pricing derivative securities. Here we recall only two such works, \cite{Guo} and \cite{Eli}. In \cite{Guo}, the author considers a stock price model which allows for the drift and the volatility coefficients to switch according to two-states. This market is incomplete, but is completed with
 new securities. In \cite{Eli} the problem of option pricing is considered in a model where the risky underlying assets are driven by Markov-modulated Geometric Brownian motions. A regime switching Esscher transform is used in order to find a martingale pricing measure. When it comes to optimal investment in regime switching
 markets we point to \cite{Cad}. In their paper they allow for the risk preference to switch according
 to the regime.

 The discount rate in our paper is stochastic, exogenous and it depends on the regime. By the best of
 our knowledge is the first work to consider stochastic discounting within the Merton problem framework.
 In a discrete time model, \cite{Tosh} considers a cyclical discount factor.

 Non constant discount rates lead to time inconsistency of the decision maker as shown in \cite{EkePir}
 and \cite{EkePirMbo}. The resolution is to consider subgame perfect equilibrium strategies. These are
 strategies which are optimal to implement now given that they will be implemented in the future. After
 we introduce this concept we try to characterize it. In order to achieve this goal the methodology
 developed in  \cite{EkePirMbo} is employed. That is a new result in stochastic control theory: it mixes the
 idea of value function (from the dynamic programming principle) with the idea that in the future ``optimal
 trading strategies'' are implemented (from the maximum principle of Pontryagin). The new twist in our paper is
 the Markov chain, and the mathematical ingredient used is $It\hat{o}'s$ formula for the Markov-modulated diffusions. Thus, we obtain a system  of four equations: first equation says that the value function is equal to the continuation utility of subgame perfect strategies; second equation is the wealth equation generated by subgame perfect strategies; the last two equations relate the value function to the subgame
 perfect strategies. The end result is a complicated system of PDEs, SDE and a nonlinear equation with
 a nonlocal term. The investor's  risk preference in this model is of CRRA type which suggest an ansatz for
 the value function (it disentangles the time, the space and the Markov chain state component). This result
 in subgame perfect strategies which are time/state dependent and linear in wealth. In the special case
 of logarithmic utility we can compute them explicitly. If constant discount rates, we notice that
  subgame perfect strategies coincide with the optimal ones.

 The goal of this paper is twofold: first, to consider a model with stochastic discount rates and second, to
 study the relationship between consumption and discount rates (which resembles IMI and/or DMI ). In
 the Merton problem with constant discount rates we show that higher the discount rate higher the consumption
 rate (this is somehow an inverse relationship to IMI). We explore this relationship in our model with stochastic
 discount rates. Numerical experiments revealed that the consumption rate is higher in the market states with
 higher discount rate. We provide an analytic proof of this result for the special case of logarithmic utility. This result is somehow consistent with the discount rate monotonocity of consumption rate in the Merton problem. It can also explain the IMI effect: if we observe the consumption rate, then a possible upside jump can be linked to a jump in the discount rate and vice versa. The effect of risk aversion on the consumption
 rate is analysed. Here the results are consistent with \cite{Mer69} and \cite{Mer71}. That is consumption rate is increasing
 in time for most levels of risk aversion except at very high levels when is decreasing ( this can be explained by
 the investor's increased appetite for risk which leads to more investment in the risky asset and a deceasing
 consumption rate).

 The reminder of this paper is organized as follow.
 In section $2$ we describe the model and formulate the objective.
 Section $3$ contains the main result under the power utility and logarithmic utility.
 In section 4, we present the numerical results. Section 5 examines consumption versus discount rate.
  The paper ends with an appendix containing the proofs.

\section{The Model}

\subsection{The Financial Market}
Consider a probability space $(\Omega,\{\mathcal{F}_t\}_{0\leq t\leq
T},\mathcal{F},\P),$ which accommodates a standard Brownian
motion $W=\{W(t),{t\geq0}\}$ and a homogeneous finite state
continuous time Markov Chain (MC) $J=\{J(t), t\geq0\}.$ For simplicity assume that MC takes
values in $\mathcal{S}=\{\textbf{0},\textbf{1}\}.$ Our results hold true in the more general situation
of  $\mathcal{S}$ having finitely many states. The filtration $\{\mathcal{F}_t\}_{0\leq
t\leq \infty}$ is the completed
 filtration generated by $\{W(t)\}_{t\in[0,\infty)}$ and
 $\{J(t)\}_{t\in[0,\infty)}, $ that is $\mathcal{F}_t=\mathcal{F}_t^J
 \bigvee \mathcal{F}_t^W.$
  We assume that the stochastic processes
$W$ and $J$ are independent. The MC $J$ has a generator
$\Lambda=[\lambda_{ij}]_{\mathcal{S}\times\mathcal{S}}$ with
$\lambda_{ij}\geq0$\ for $i\neq j,$ and  $\sum_{j\in
\mathcal{S}}\lambda_{ij}=0$ for every $i\in \mathcal{S}.$

In our setup the financial market consists of a bank account $B$ and a
risky asset $S$, that are traded continuously over a finite time
horizon $[0,T]$ (here $T\in (0,\infty)$ is an exogenously given deterministic time).
 The price process of the bank account and risky asset are governed
 by the following Markov-modulated SDE:
 \begin{eqnarray}
 &&dB(t)=r(t,{J(t)})B(t)dt,\nonumber\\
&&dS(t)=S(t)\left[\alpha(t,{J(t)})\,dt
+\sigma(t,{J(t)})\,dW(t)\right],\quad0\leq t\leq \infty,\nonumber
 \end{eqnarray}
 where $B(0)=1$ and $S(0)=s>0$ are the initial prices. The functions
 $r(t, i), \alpha(t,i), \sigma(t,i):  i\in \mathcal{S},$ are assumed to be
 deterministic, positive and continuous in $t.$ Given the state $i$ of
 the MC at $t$ they represent the riskless rate, the stock volatility and the stock return.
 Moreover
 $$\mu(t,i)\triangleq\alpha(t,i)-r(t,i) $$
 stands for the stock excess return.

\subsection{Investment-consumption strategies and wealth processes}
In our model, a representative investor continuously invests in the stock,bond and consumes. An acceptable
investment-consumption strategy is defined below:
\begin{definition}
\label{def:portfolio-proportions} An $\R^{2}$-valued stochastic
process $\{u(t):=(\pi(t), c(t))\}_{t\in[0,\infty)}$ is called an
  admissible strategy process and write $u\in \mathcal{A}$ if
it is ${\mathcal{F}_t}-$ progressively measurable and it satisfies the following integrability condition
\begin{equation}%
\label{kj***} E[\int_0^t |\pi(s)\mu(s,J(s))-c(s)|\, ds+\int_0^t
      |\pi(s) \sigma({s,J(s)})|^2\, ds]<\infty, \text{ a.s., for all }
 t\in [0,\infty).
\end{equation}
\end{definition}
Here $\pi(t)$ stands for the dollar value invested in stock at time $t$ and
$c(t)\geq0$ for the consumption. ${X^{u}(t)}:=X(t)$ represents the
wealth of the investor at time $t$ associated with the trading strategy $u;$
it satisfies the following stochastic differential equation (SDE)
\begin{eqnarray}\label{equ:wealth-one}
dX(t)&=& \left(r(t,{J(t)}) X^{u}(t)+\mu(t,{J(t)})\pi(t)-c(t))\,dt
+\sigma(t,{J(t)})\pi(t)\,dW(t),\right.
\end{eqnarray}
where $X(0)=x>0$ is the initial wealth and $J_0=i\in{\mathcal{S}}$ is initial state.
  This SDE is called the self-financing condition. Under the regularity condition \eqref{kj***} imposed on $\{\pi(t),
c(t)\}_{t\in[0,\infty)}$ above, the SDE
 \eqref{equ:wealth-one} admits a unique strong solution. In the end of this section, we introduce
 further assumptions.
 \begin{assumption}  The utility function of
the investor is of CRRA type, i.e.,
$$U(x)=\frac{x^\gamma}{\gamma}$$where $\gamma<1, \gamma\neq0.$
\end{assumption}

\noindent Thus, the inverse marginal utility function is
\begin{equation}\label{mar}
I(x)\triangleq(U')^{-1}(x)=x^{\frac{1}{\gamma-1}}.
\end{equation}

\begin{assumption}
For any risk aversion level $\gamma<1,$ the following inequalities hold
\begin{equation}\label{089}
\mathbb{E} \sup_{t\in[0,T]}|X(t)|^\gamma<\infty.
\end{equation}
\begin{equation}\label{189}
\mathbb{E} \sup_{t\in[0,T]}|c(t)|^\gamma<\infty.
\end{equation}
\end{assumption}

\subsection{The discount rate}
As we mentioned in the introduction, this paper considers stochastic discount rates.
An easy way to achieve this is to let the discount rate depend on the state of MC. Thus,
 at some intermediate time $t\in [0, T]$ the discount are is $\rho_{J(t)},$ for some positive
 constants $\rho_{\bf{0}}$ and $\rho_{\bf{1}}.$ The intuition of this way of modeling discount
 rate stems from the connection between market state and discount rates (this can be explained
 by some models with endogenous discount rate which may be influenced by economic factors).

\subsection{The Risk Criterion}
In our model, the investor decides what investment/consumption
strategy to choose according to the expected utility risk criterion.
Thus, investor's goal is to maximize utility of intertemporal consumption
and final wealth. The novelty here is that we allow investor to update
the risk criterion and to reconsider the optimal strategies she/he computed in the past.
This will lead to a time inconsistent behavior as we show below. Let the agent start with a given
positive wealth $x,$ and a given market state $i,$ at some instant
$t.$ The optimal $t-$trading strategy
 $\{\tilde{\pi}_{t}({s}),\tilde{c}_{t}({s})\}_{s\in[t,T]}$ is chosen
such that
\begin{eqnarray}
\sup_{u\in\mathcal{A}}\mathbb{E}\left[\int_{t}^{T}e^{-\rho_i(u-t)}
U(c(u))\,du
+e^{-\rho_i(T-t)}{U}(X(T))|X(t)=x,J(t)=i\right]\nonumber\\
=\mathbb{E}_t^{x,i}\left[\int_{t}^{T}e^{-\rho_i(u-t)}U(\tilde{c}_{t}(u))\,du
+e^{-\rho_i(T-t)}{U}(\tilde{X}(T))\right]\nonumber.
\end{eqnarray}
Throughout the paper we denote $\mathbb{E}_t^{x,i}\left[ \right] \triangleq \mathbb{E}\left[ |X(t)=x,J(t)=i\right].$
 The optimal $t-$trading strategy
 $\{\tilde{\pi}_{t}({s}),\tilde{c}_{t}({s})\}_{s\in[t,T]}$ is derived by the supermartingale/martingale principle. This
 leads to the Hamilton-Jacobi-Bellman (HJB) equation
\begin{equation}\label{hjb}
\frac{\partial V}{\partial
s}(t,s,x,i)+\sup_{\pi,c}\left[(r x+\mu \pi-c) \frac{\partial
V}{\partial x}(t,s,x,i)+\frac{1}{2}\sigma^2\pi^{2}
\frac{\partial^{2} V}{\partial x^{2}}(t,s,x,i)+U(c)\right]\end{equation}$$-
\rho_iV(t,s,x,i)+\sum_{j\in \mathcal{S}}\lambda_{ij}V(t,s,x,j)=0,
$$
with the boundary condition
\begin{equation}\label{boundarycondition}
V(t,T,x,i)={U}(x).\ \ \ 
\end{equation}
Here $i$ stands for the value of MC at time $t.$ Thus, the HJB depends on the current time $t$
( through $\rho_i$) and this dependence is inherited by the $t-$optimal trading strategy. This in turn
leads to time inconsistencies. The resolution is to introduce subgame perfect equilibrium strategies. They are optimal now given that they will be implemented in the future.

\section{The Main Result}

\subsection{The subgame perfect trading strategies}
For a policy process $\{u(t)\triangleq{\pi}(t),{c}(t)\}_{t\in[0,T]}$
 satisfying \eqref{kj***} and its corresponding wealth process
$\{X(t)\}_{t\in[0,T]}$ given by  \eqref{equ:wealth-one}, we denote
the expected utility functional by
\begin{equation}\label{01FUNCT}
J(t,x,i,\pi,c)\triangleq\mathbb{E}_t^{x,i}\left[\int_{t}^{T}e^{-\rho_i(s-t)}
U(c(s)) \,ds+e^{-\rho_i(T-t)}{U}(X(T))\right].
\end{equation}
 Following \cite{EkePir} we shall give a rigorous mathematical
formulation of the equilibrium strategies in the formal definition
below.

\begin{definition}\label{finiteh}
Let $F=(F_{1},F_{2}):[0,T]\times R^+\times\mathcal{S}\rightarrow
{R}^+\times\mathcal{S}$ be a map such that for any $t,x>0$ and $i\in
\mathcal{S}$
\begin{equation}\label{opt}
{\lim\inf_{\epsilon\downarrow 0}}\frac{J(t,x,i,F_{1},F_{2})-
J(t,x,i,\pi_{\epsilon},c_{\epsilon})}{\epsilon}\geq 0,
\end{equation}
where
$$J(t,x,i,F_{1},F_{2})\triangleq J(t,x,i,\bar{\pi},\bar{c}),$$
\begin{equation}\label{0000eq}
\bar{\pi}(s)\triangleq{F_{1}(s,\bar{X}(s),J(s))},
\quad\bar{c}(s)\triangleq{F_{2}(s,\bar{X}(s),J(s))},
\end{equation}
and $\{\bar{\pi}(s),\bar{c}(s)\}_{s\in[t,T]}$  satisfies
\eqref{kj***}.
 Here, the  process $\{\bar{X}(s)\}_{s\in[t,T]}$ is the wealth corresponding to  $\{\bar{\pi}(s),\bar{c}(s)\}_{s\in[t,T]}.$
The process $\{{\pi}_{\epsilon}(s),{c}_{\epsilon}(s)\}_{s\in[t,T]}$
 is another investment-consumption strategy defined by
\begin{equation}\label{1e}
\pi_{\epsilon}(s)=\begin{cases} \bar{\pi}(s),\quad
s\in[t,T]\backslash E_{\epsilon,t}\\
\pi(s), \quad s\in E_{\epsilon,t}, \end{cases}
\end{equation}

\begin{equation}\label{2e}
c_{\epsilon}(s)=\begin{cases} \bar{c}(s),\quad s\in[t,T]\backslash E_{\epsilon,t}\\
c(s), \quad s\in E_{\epsilon,t}, \end{cases}
\end{equation}
with $E_{\epsilon,t}=[t,t+\epsilon];$ $\{{\pi}(s),{c}(s)\}_{s\in
E_{\epsilon,t} }$ is any trading strategy for which
$\{{\pi}_{\epsilon}(s),{c}_{\epsilon} (s)\}_{s\in[t,T]}$ is an
admissible policy. If \eqref{opt} holds true, then  $\{\bar{\pi}(s),\bar{c}(s)\}_{s\in[t,T]}$ is a subgame perfect strategy.
\end{definition}

\subsection{The value function}
Our goal is in a first step to characterize the subgame perfect strategies and then to find them.
Inspired by \cite{EkePir}, the value function $v$ satisfies
\begin{equation}\label{00ie1}
v(t,x,i)=\mathbb{E}_t^{x,i}\left[\int_{t}^{T}{e^{-\rho_i(s-t)}}
U(F_{2}(s,\bar{X}(s),J(s)))\,ds+
{e^{-\rho_i(T-t)}}{U}(\bar{X}(T))\right].
\end{equation}
 Recall that $\{\bar{X}(s)\}_{s\in[0,T]}$
is the wealth process  corresponding to  $\{\bar{\pi}(s),\bar{c}(s)\}_{s\in[t,T]},$ so it solves the SDE
\begin{equation}\label{*10dyn}
\!\!\!\!d\bar{X}(s)\!\!=\!\![r(s,J(s))\bar{X}(s)\!+\!\mu(s,J(s))
F_{1}(s,\bar{X}(s),J(s))\!-\!F_{2}(s,\bar{X}(s),J(s))]ds\!+\! \sigma(s,J(s))
F_{1}(s,\bar{X}(s))dW(s).
\end{equation}
 Moreover, $F=(F_{1},F_{2})$ is given by
 \begin{equation}\label{109con}
 F_{1}(t,x,i)=-\frac{\mu(t,i)\frac{\partial v}{\partial x}(t,x,i)}
 {\sigma^2(t,i)
\frac{\partial^{2} v}{\partial x^{2}}(t,x,i)},\,\,
F_{2}(t,x,i)=I\left(\frac{\partial v}{\partial
x}(t,x,i)\right),\,\,\,t\in[0,T].
\end{equation}
Thus, the value function is characterized by a system of four equations: one integral equation
with nonlocal term \eqref{00ie1}, one SDE \eqref{*10dyn} and two PDEs \eqref{109con}. Of course the existence of such a function $v$ satisfying the
equations above is not a trivial issue. We will take
advantage of the special form of the utility function to simplify the
problem of finding $v.$ We search for: \begin{equation}\label{oP22}
  v(t,x,i)=g(t,i)\frac{x^\gamma}{\gamma},\ \ x\geq0,
\end{equation}
where the function $g(t,i)$ is to be found. We consider the
case $\gamma\neq 0$ (the case of logarithmic utility will be treated
separately). In the light of equations (\ref{109con}) one gets
\begin{equation}\label{eQ1}
 F_1(t,x,i)=\frac{\mu(t,i) x}{\sigma^2(t,i) (1-\gamma)},\ \
 F_2(t,x,i)={g^{\frac{1}{\gamma-1}}(t,i)}x.
\end{equation}
By (\ref{*10dyn}), the associated wealth process satisfies  the
following SDE:
\begin{eqnarray}\label{dyn1}
d\bar{X}(s)&=&\left[r(s,J(s))+\frac{\mu^2(s,J(s))}{\sigma^2(s,J(s))(1-\gamma)}
-g^{\frac{1}{\gamma-1}}(s,J(s))\right]\bar{X}(s)ds\nonumber\\
&&+\frac{\mu(s,J(s))}{\sigma(s,J(s))(1-\gamma)}\bar{X}(s)dW(s).
\end{eqnarray}

This is a linear SDE which can be easily solved. By plugging $v$ of (\ref{oP22}) into
(\ref{00ie1}) (with $F_1,F_2$ of (\ref{eQ1}) and $\bar{X}$ of (\ref{dyn1}), we obtain the following equation for
$g(t,i), i\in \mathcal{S}:$
\begin{equation}\label{oDe} \frac{\partial g}{\partial t}(t,i)+[\gamma r(t,i)
+\frac{\mu^2(t,i)\gamma}{2\sigma^2(t,i)(1-\gamma)}-\rho_i]g(t,i)+
\sum_{j\in
S}\lambda_{ij}g(t,j)+(1-\gamma)g^{\frac{\gamma}{\gamma-1}}(t,i)=0,
\end{equation} with the
final condition $g(T,i)=1.$ Next we show that there exists a unique solution of this ODE system. Let us summarize these findings:
 \begin{lemma}\label{l4.1} There exists a unique continuously differentiable
 solution $g(t,i), i\in \mathcal{S}$ of the system (\ref{oDe}). Furthermore, $v(t,x,i)=g(t,i)\frac{x^\gamma}{\gamma}$ is a value function, meaning that $v$ solves (\ref{00ie1}) with $F_1, F_2$ of
(\ref{eQ1}) and $\bar{X}$ of (\ref{dyn1}).\end{lemma}

Appendix A proves this Lemma.
\begin{flushright}
$\square$
\end{flushright}


The following Theorem states the central result of our paper.

\begin{theorem}\label{4.1}
 Suppose that $v(t, x, i)$ is given by (\ref{oP22}) with $g(t,i), i\in \mathcal{S}$ the solution of the system (\ref{oDe}). Let $\bar{X}$ be the solution of SDE \eqref{dyn1}. Then $\{\bar{\pi}(t),\bar{c}(t)\}_{s\in[0,T]}$ given by
 \begin{equation}\label{eQ}
 \bar{\pi}(t)=\frac{\mu(t,J(t)) X(t)}{\sigma^2(t,J(t)) (1-\gamma)},\ \
\bar{c}(t)={g^{\frac{1}{\gamma-1}}(t,J(t))} {X(t)},
\end{equation}
is a subgame perfect strategy.
\end{theorem}
Appendix B proves this Theorem.
\begin{flushright}
$\square$
\end{flushright}

\begin{remark}
In the case of constant discount rate, i.e., $\rho_{\bf{0}}=\rho_{\bf{1}},$ the subgame
perfect strategies coincide with optimal ones. This can be seen by looking at the
integral equation \eqref{zhang} which is exactly HJB \eqref{hjb}
(after the first order conditions are implemented).

\end{remark}

Next we turn to the special case of logarithmic utility.

\subsection{Logarithmic Utility}
When the risk aversion $\gamma=0,$ we search for the value function $v$
 of the following form:
\begin{equation}\label{oP}
  v(t,x,i)=h(t,i)\log(x)+l(t,i).
  \end{equation}
Arguing as in the previous subsection we find that
the functions $h,\
l:[0,T]\times\mathcal{S}\rightarrow R^{+}$ should satisfy the following system of
equations:
 \begin{equation}\label{oDeh} \frac{\partial h}{\partial t}(t,i)-\rho_i h(t,i)+
\sum_{j\in S}\lambda_{ij}h(t,j)+1=0,
\end{equation}
$$\frac{\partial l}{\partial t}(t,i)+(r(t,i)+\frac{\mu^2(t,i)}{2\sigma^2(t,i)})h(t,i)-\log l(t,i)-\rho_i l(t,i)
+\sum_{j\in S}\lambda_{ij}l(t,j)-1=0,$$
 with the
final conditions $h(T,i)=1$ and $l(T,i)=0.$ Notice that $h$ solves a linear ODE system and it can be
found explicitly. With $h$ known, $l$ also solves  a linear ODE system. The subgame perfect strategy $\{\bar{\pi}(t),\bar{c}(t)\}_{s\in[0,T]}$ is given by
 \begin{equation}\label{eQ}
 \bar{\pi}(t)=F_{1}(t,J(t)) X(t)),\qquad
\bar{c}(t)=F_{2}(t,J(t)) X(t)),
\end{equation}
with
\begin{equation}\label{eQ10}
 F_1(t,x,i)=\frac{\mu(t,i) x}{\sigma^2(t,i) },\ \
 F_2(t,x,i)=h^{-1}(t,i)x.
\end{equation}

Next we want to explore the relationship between subgame perfect consumption and discount rate. The following lemma shows that the higher discount rate the
more consumption.

\begin{lemma}\label{4.4}
 Assume that $\rho_{\bf{0}}>\rho_{\bf{1}},$ and $h:[0,T]\times\mathcal{S}\rightarrow R^{+}$  solves
 the linear ODE (\ref{oDeh}). Then $F_2(t,x,{\bf{0}})>F_2(t,x,{\bf{1}})$ for any $x>0.$
 In other words the subgame perfect consumption rate is higher in the states in which the discount rate is higher. \end{lemma}
Appendix C proves this Theorem.

\begin{flushright}
$\square$
\end{flushright}

\section{Numerical Analysis }

In this section, we use Matlab's powerful ODE solvers (especially the functions  ode23 and ode45) to perform
numerical experiments. We numerically solve ODE system \eqref{oDe} and this in turn yield the subgame
perfect strategies. Let the market coefficients be $\mu_1=0.15,\ \mu_2=0.15,\
\sigma_1=0.25,\ \sigma_2=0.25,\ \ r_1=0.05,\ r_2=0.05;\ $ $
\gamma=0.7,0,-0.5$ and $-1,$  the discount rate $ \rho_1=0.9,\
\rho_2=0.3.$ We take the Markov Chain generator to be
\begin{equation}
\left(
\begin{array}{cc}
 -6.04 & 6.04 \\
 10.9 & -10.9
 \end{array}
\right)\nonumber
\end{equation}
 Next, define the consumption rate by
$$C(t,J(t)):=\frac{F_2(t,X(t),J(t))}{X(t)}=g^{\frac{1}{\gamma-1}}(t,J(t)),$$
with the functions $g(t,i)$ for $i=1,2$ being the solution of (\ref{oDe}).

\begin{center}

\begin{tabular}{cc}

\includegraphics
[scale=.95]{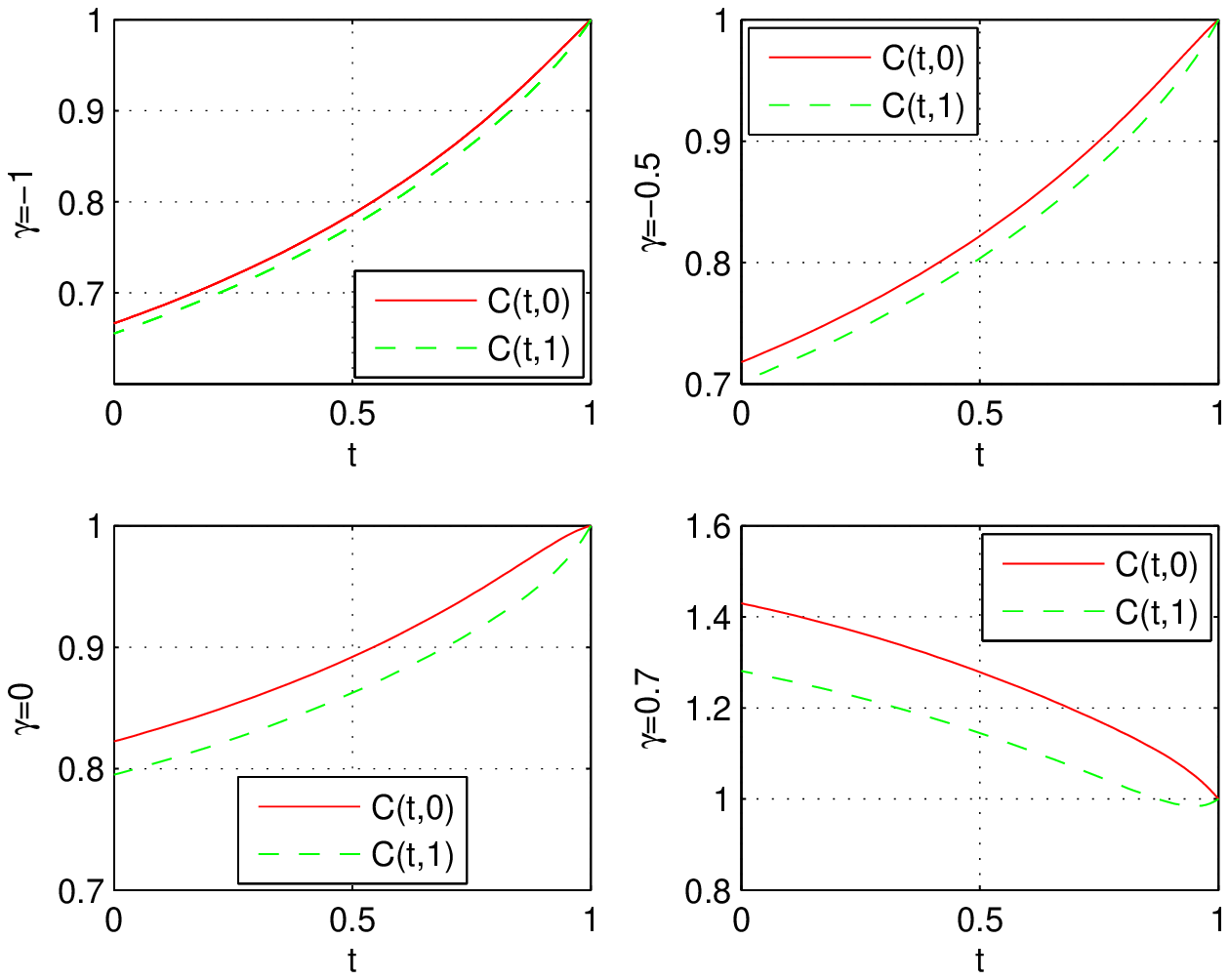}
\end{tabular}

\end{center}
$$\mbox{Fig.}\,\,1 $$

\noindent Fig. $1.$ Equilibrium proportion of wealth consumed for
different $\gamma$ and $\rho.$ The $x$ axis represents the time, and the $y$ axis
 the consumption rate $C.$

 \begin{remark}
 As $t\rightarrow T$ consumption rate approaches $1$ as it was to be
expected. From the picture with $\gamma=0.7$, we can see that
the consumption rate decreases with time. This is explained by the fact that the
higher $\gamma,$ the less risk aversion, which leads to higher
proportion of wealth invested into stock and less consumption. This is consistent with a graph from \cite{Mer69}. We see from the graphs that consumption rate increases when $\gamma$ increases
(this is also consistent with graphs from \cite{Mer69}). For a given risk aversion level
$\gamma$, a higher discount rate will result in higher consumption rate. The difference in the consumption
rates in the two MC states is decreasing with respect to $\gamma.$

\end{remark}

\section{Consumption versus Discount Rate}
We were able to prove for the case of logarithmic utility
that higher discount rates lead to higher consumption
rates. Numerical evidence suggest that this is also true
for general power utilities. On the other hand, in a model with
constant discount rate the same result holds true as the following
Lemma shows.

\begin{lemma}\label{merton}
Let $\bar{C}(t), t\in[0,T]$ be the optimal consumption rate in a model  with
constant discount rate $\rho.$ Then
$$\frac{\partial C(t)}{\partial \rho}>0, \qquad t\in [0,T].$$
\end{lemma}

Appendix D proves this Theorem.

\begin{flushright}
$\square$
\end{flushright}

\section{Appendix}

\subsection{It$\hat{o}$'s formula for Markov Chain modulated diffusions}

 Suppose the stochastic processes $X(t)$ satisfies the SDE
$$dX(t)=\mu(t,X(t),J(t))dt+\sigma(t,X(t),J(t))dW(t)$$
$$X(0)=x, a.s.,$$
for some $x\in R,$ and $G(\cdot,\cdot,i)\in C^{1,3}([0,T]\times R)$
for each $i\in \mathcal{S}.$ Then
\begin{eqnarray}
G(t,X(t),J(t))&=&G(0,x,J(0))+\int_0^t \Gamma
G(s,X(s),J(s))ds\nonumber\\&&+\frac{\partial G}{\partial
x}(s,X(s),J(s))\sigma(s,X(s),J(s))dW(s)\nonumber\\&&+\sum_{J(t)\neq
i}\int_0^t (G(s,X(s),J(t))-G(s,X(s),i))dM(t),\nonumber
\end{eqnarray}
where
\begin{eqnarray}
\Gamma G(t,x,i)&=&\frac{\partial G}{\partial
t}(t,x,i)+\mu(t,x,i)\frac{\partial G}{\partial
x}(t,x,i)+\frac{1}{2}\sigma^{2}(t,x,i)\frac{\partial^2G}{\partial
x^2}+\sum_{j\in \mathcal{S}}\lambda_{ij}
(G(t,x,j)-G(t,x,i)).\nonumber
\end{eqnarray}
Here $M(t),\,\, t\in[0,T]$ is the martingale process associated with the Markov Chain.

\subsection{A Proof of Lemma \ref{l4.1} }
The existence of a unique solution $g$ of ODE system (\ref{oDe}) is granted
locally in time by a fixed point theorem. If we can establish estimates for $g$ then
this local solution is also global solution. Let us introduce the process $\{M(v),v\geqslant t\}$ by
\begin{equation}\label{M}
M(v)\triangleq g(v,J_v)-g(t,J_t)-\int_t^v(g_u(u,J_u)+\sum_{j\in\mathcal{S}}\lambda
_{J_uj}g(u,j))du,
\end{equation}
 with $g$ of  (\ref{oDe}). Dynkin formula implies that the process
$\{M(v),v\geqslant t\}$ defined by (\ref{M}) is a martingale. Further let
\begin{equation}\label{K}
K(v)\triangleq\exp\left\{\int_t^v \left[\gamma r(u,J(u))
+\frac{\mu^2(u,J(u))\gamma}{2\sigma^2(u,J(u))(1-\gamma)}-\rho_{J(u)} \right]du\right\}.
\end{equation}
By product rule
\begin{eqnarray}\label{Kg}
d(K(v)g(v,J_v))
&=&K(v)\left[d g(v,J_v)+g(v,J_v)\left(\gamma r(v,J_v)
+\frac{\mu^2(v,J_v)\gamma}{2\sigma^2(v,J_v)(1-\gamma)}-\rho_{J_v} \right)dv\right]\nonumber\\
&=&K(v)[d M(v)+(\gamma-1)g^{\frac{\gamma}{\gamma-1}}(v,J_v)dv].
\end{eqnarray}
Integrating (\ref{Kg}) from $t$ to $T$, we get
\begin{equation}\label{KGI}
K(T)-g(t,J_t)=\int_t^TK(v)dM(v)+(\gamma-1)\int_t^T
K(v)g^{\frac{\gamma}{\gamma-1}}(v,J_v)dv.\end{equation} Taking
expectation on the both sides of (\ref{KGI}) and letting $J_t=i$, it leads to
\begin{equation}\label{EKGI}
g(t,i)=\mathbb{E}_{t}^iK(T)+(1-\gamma)\mathbb{E}_{t}^i \int_t^T
K(v)g^{\frac{\gamma}{\gamma-1}}(v,J_v)dv.\end{equation}
Boundedness assumption on market coefficients makes the process $\{K(v),v\geqslant t\}$
of \eqref{Kg} bounded. This in turn yields that $g(t,i)$ is uniformly bounded
from below. Next we want to prove that is also bounded from above. For a vector $y=(y_1, y_2,...,y_n)$
in $\mathbf{R}^{n}$ we introduce the $\parallel \cdot \parallel_1$ norm by :
 $$\parallel y \parallel_1=\sum_{i=1}^n|y_i|.$$
When $\gamma$ is positive \eqref{EKGI} yields an upper
bound on $g$ (since $g$ is bounded from below). When $\gamma<0,$ because $0<\frac{\gamma}{\gamma-1}<1$ we get from  \eqref{EKGI} that

\begin{equation}
g(t,0)<\mathbb{E}_{t}^0K(T)+\mathbb{E}_{t}^0(1-\gamma)\int_t^T
K(v)g(v,J_v)dv.\end{equation}

\begin{equation}
g(t,1)<\mathbb{E}_{t}^1K(T)+\mathbb{E}_{t}^1(1-\gamma)\int_t^T
K(v)g(v,J_v)dv.\end{equation} Thus,
$$||g(t,\cdot)||_1\leq 2\bar{K}+(1-\gamma)\bar{K}\int_t^T||g(v,\cdot)||_1dv,$$
for some positive constant $\bar{K}.$ Gronwal's inequality yields an upper bound on $g.$
Next we prove uniqueness of the solution of (\ref{oDe}).
Indeed, let $g_1$ and $g_2$ be two solutions of (\ref{oDe}), then by (\ref{EKGI}) it follows that
\begin{eqnarray}\label{unique}
|g_1(t,i)-g_2(t,i)|&=&(1-\gamma)|\mathbb{E}_t^i\int_t^T
K(v)g_1^{\frac{\gamma}{\gamma-1}}(v,J_v)dv -\mathbb{E}_t^i\int_t^T
K(v)g_2^{\frac{\gamma}{\gamma-1}}(v,J_v)dv|\nonumber\\
&&\leq \mathbb{E}_t^i\int_t^T K(v)|g_1^{\frac{\gamma}{\gamma-1}}
(v,J_v)-g_2^{\frac{\gamma}{\gamma-1}}(v,J_v)|dv\nonumber\\
&&\leq C_1 \mathbb{E}_t^i\int_t^T K(v)|g_1
(v,J_v)-g_2(v,J_v)|dv\nonumber\\
&&\leq C_2\mathbb{E}_t^i\int_t^T |g_1
(v,J_v)-g_2(v,J_v)|dv\nonumber\\
&&\leq C_2 \int_t^T (|g_1 (v,i)-g_2(v,i)|+|g_1
(v,j)-g_2(v,j)|)dv\nonumber\\
&&= C_2 \int_t^T ||g_1 (v,\cdot)-g_2(v,\cdot)||_1   dv,
\end{eqnarray}
for some positive constants $C_1$ and $C_2.$ Consequently

$$||g_1(t,\cdot)-g_2(t,\cdot)||\leq C_2 \int_t^T ||g_1
(v,\cdot)-g_2(v,\cdot)||dv.$$Thus, by Gronwal's inequality it follows that
$$g_1(t,\cdot)=g_2(t,\cdot).$$
Next we want to prove that $v(t,x,i)=g(t,i)\frac{x^\gamma}{\gamma}$ is a value function, i.e., $v$ solves (\ref{00ie1}) with $F_1, F_2$ of
(\ref{eQ1}) and $\bar{X}$ of (\ref{dyn1}). Since the function
$v$ is sufficiently differentiable, we can compute the derivative of
$v$:
$$v_t(t,x,i)=g_t(t,i)\frac{x^\gamma}{\gamma},\,\,v_x(t,x,i)=g(t,i){x^{\gamma-1}},\,\,
v_{xx}(t,x,i)=(\gamma-1)g(t,i){x^{\gamma-2}}.$$ Next we show that $v$ solves the PDE system \eqref{zhang}
(which Lemma \ref{main} shows that is equivalent to \eqref{00ie1}). Substituting the above derivatives  into the equation (\ref{zhang})
yields
\begin{eqnarray}\label{traian1}
&&\frac{\partial g}{\partial t}(t,i)\frac{x^\gamma}{\gamma}+[(r(t,i) x+\mu(t,i)
F_{1}(t,x,i)-F_{2}(t,x,i))g(t,i){x^{\gamma-1}} +\frac{\sigma_i^{2}
F_{1}^{2}(t,x,i)}{2} (\gamma-1)g(t,i){x^{\gamma-2}}\nonumber\\
&&+U(F_{2}(t,x,i))]+\sum_{j\in
\mathcal{S}}\lambda_{ij}g(t,i)\frac{x^\gamma}{\gamma}-\rho_i
g(t,i)\frac{x^\gamma}{\gamma}=0.\end{eqnarray}
After cancelation of $x^{\gamma}$ we recover the ODE system  (\ref{oDe}); hence $v$ solves \eqref{zhang}
(and also  \eqref{00ie1}).
\begin{flushright}
$\square$
\end{flushright}

\subsection{B Proof of Theorem \ref{4.1}}

We need the following Lemma which gives a PDE version of \eqref{00ie1}.

\begin{lemma}\label{main}For a given $i\in\mathcal{S},$
assume there exists a function $v: [0,T]\times
R^+\times\mathcal{S}\rightarrow R$ of class $C^{1,2}$ which
satisfies \eqref{00ie1}. Then $v$ solves the following equation
$$\frac{\partial v}{\partial t}(t,x,i)+(r(t,i) x+\mu(t,i) F_{1}(t,x,i)-F_{2}(t,x,i))
\frac{\partial v}{\partial x}(t,x,i) +\frac{\sigma^2(t,i)
F_{1}^{2}(t,x,i)}{2} \frac{\partial^2 v}{\partial x^2}(t,x,i)$$
\begin{equation}\label{zhang}+U(F_{2}(t,x,i))+\sum_{j\in \mathcal{S}}\lambda_{ij}v(t,x,j)-\rho_i
v(t,x,i)=0, \end{equation} with the boundary condition
$v(T,x,i)=U(x).$ 
\end{lemma}

\noindent Proof: We rewrite equation \eqref{00ie1} as
\begin{equation}\label{1ie1}
v(t,x,i)=\int_{t}^{T}e^{-\rho_i(s-t)} f(t,s,x,i)\,ds +
e^{-\rho_i(T-t)} h(t,x,i),
\end{equation}
where
$$f(t,s,x,i)\triangleq\mathbb{E}_t^{x,i}[U(F_{2}(s,\bar{X}(s),J(s)))],\quad h(t,x,i)\triangleq\mathbb{E}_t^{x,i}[U(\bar{X}(T))].$$
For a fixed time $s,$ the process $\{f(t,s,\bar{X}(t),J(t))\}_{0\leq t\leq
s}$ is a martingale. Thus, in the light of (\ref{*10dyn}) the function
$f, h$ satisfy the following PDEs :
\begin{eqnarray}\label{PDEf}
&&\frac{\partial f}{\partial t}(t,s,x,i)+(r(t,i) x+\mu(t,i)
F_{1}(t,x,i)-F_{2}(t,x,i))\frac{\partial f}{\partial
x}(t,s,x,i)\nonumber\\&+&\frac{\sigma^2(t,i) F_{1}^{2}(t,x,i)}{2}
\frac{\partial^2 f}{\partial x^2}(t,s,x,i)+\sum_{j\in
\mathcal{S}}\lambda_{ij}f(t,x,j)=0, f(t,t,x,i)=U(F_2(t,x,i))
\end{eqnarray}

\begin{eqnarray}\label{PDE1}
&&\frac{\partial h}{\partial t}(t,s,x,i)+(r(t,i) x+\mu(t,i)
F_{1}(t,x,i)-F_{2}(t,x,i))\frac{\partial h}{\partial
x}(t,s,x,i)\nonumber\\&+&\frac{\sigma^2(t,i) F_{1}^{2}(t,x,i)}{2}
\frac{\partial^2 h}{\partial x^2}(t,s,x,i)+\sum_{j\in
\mathcal{S}}\lambda_{ij}h(t,x,j)=0,\,\,h(T,x,i)=U(x).
\end{eqnarray}
By differentiating  \eqref{1ie1} with respect to $t$ we get

\begin{equation}\label{t1ie1}
\frac{\partial v}{\partial t}(t,x,i)=\int_{t}^{T}e^{-\rho_i(s-t)}
\frac{\partial f}{\partial t}(t,s,x,i)\,ds+ e^{-\rho_i(T-t)}
\frac{\partial h}{\partial t}(t,x,i)+\rho_i v(t,x,i)-f(t,t,x,i).
\end{equation}
Moreover
\begin{equation}\label{xie1}
\frac{\partial v}{\partial x}(t,x,i)=\int_{t}^{T}e^{-\rho_i(s-t)}
\frac{\partial f}{\partial x}(t,s,x,i)\,ds
+e^{-\rho_i(T-t)}\frac{\partial h}{\partial x}(t,x,i).
\end{equation}

\begin{equation}\label{yie1}
\frac{\partial^2 v}{\partial x^2}(t,x)=\int_{t}^{T}e^{-\rho_i(s-t)}
\frac{\partial^2 f}{\partial x^2}(t,s,x,i)\,ds
+e^{-\rho_i(T-t)}\frac{\partial^2 h}{\partial x^2}(t,x,i).
\end{equation}
In the light of \eqref{PDEf}, \eqref{PDE1}, \eqref{t1ie1},
\eqref{xie1} and \eqref{yie1} it follows that

$$\frac{\partial v}{\partial t}(t,x,i)+(r(t,i) x+\mu(t,i) F_{1}(t,x,i)-F_{2}(t,x,i))
\frac{\partial v}{\partial x}(t,x,i) +\frac{\sigma^2(t,i)
F_{1}^{2}(t,x,i)}{2} \frac{\partial^2 v}{\partial x^2}(t,x,i)$$
$$+U(F_{2}(t,x,i))+\sum_{j\in
\mathcal{S}}\lambda_{ij}v(t,x,j)-\rho_i v(t,x,i)=0.
$$

\begin{flushright}
$\square$
\end{flushright}

Let us move to the proof of the Theorem  \ref{4.1}. First let us show that
 \begin{equation}\label{eQ11}
 \bar{\pi}(t)=F_{1}(t,J(t),X(t))= \frac{\mu(t,J(t)) X(t)}{\sigma^2(t,J(t)) (1-\gamma)},\ \
\bar{c}(t)=F_{2}(t,J(t)), X(t))={g^{\frac{1}{\gamma-1}}(t,J(t))} {X(t)},
\end{equation}
is admissible and it satisfies \eqref{089} and \eqref{189}. This claim follows easily
since $\bar{X}$ of \eqref{dyn1} has finite moments of any order. Next, let us define

\begin{eqnarray}
\overline{\Gamma} v(t,x,i)&\triangleq&\frac{\partial v}{\partial
t}(t,x,i)+\left(r(t,i)x-I\left(\frac{\partial v}{\partial
x}(t,x,i)\right)\right)\frac{\partial v}{\partial
x}(t,x,i)\nonumber\\&-&
\frac{\mu^{2}(t,i)}{2\sigma^{2}(t,i)}\frac{{[\frac{\partial
v}{\partial x}}(t,x,i)]^{2}}{\frac{\partial^{2} v}{\partial
x^{2}}(t,x,i)}+\sum_{i\in\mathcal{S}}\lambda_{ij}(v(t,x,i)-v(t,x,j))+U(c).
\end{eqnarray}

\begin{eqnarray}
{\Gamma}^{\pi, c} v(t,x,i)&\triangleq&\frac{\partial v}{\partial
t}(t,x,i)+\left(r(t,i)x-\mu(t,i)\pi-c\right)\frac{\partial
v}{\partial x}(t,x,i)\nonumber\\&+&\frac{1}{2}
{\sigma^{2}(t,i)}\pi^2{\frac{\partial^{2} v}{\partial
x^{2}}(t,x,i)}+\sum_{i\in\mathcal{S}}\lambda_{ij}(v(t,x,i)-v(t,x,j))+U(F_{2}(t,x,i)).
\end{eqnarray}
Then by the concavity of $v$ and the first order conditions it follows that

\begin{equation}\label{21}
\overline{\Gamma} v(t,x,i)=\max_{\{\pi,c \}} {\Gamma}^{\pi, c} v(t,x,i),\qquad (F_{1} (t,x,i), F_{2} (t,x,i))=\arg\max_{\{\pi,c \}}{\Gamma}^{\pi, c} v(t,x,i).
\end{equation}
Moreover, equation \eqref{zhang} can be written

\begin{equation}\label{31}
\overline{\Gamma} v(t,x,i)=\rho_i v(t,x,i).
\end{equation}
Let us recall that
\begin{equation}\label{yy}
J(t,x,i,F_{1},F_{2})=v(t,x,i).
\end{equation}
Thus
\begin{eqnarray}\label{7.16}
&&J(t,x,i,F_{1},F_{2})-J(t,x,i,\pi_{\epsilon},c_{\epsilon})\nonumber\\
&=&\mathbb{E}_t^{x,i}\left[\int_{t}^{t+\epsilon}e^{-\rho_i(s-t)}
[U(F_2(s,\bar{X}(s),J(s)))-U(c(s))]\,ds\right]\nonumber\\
&+&\mathbb{E}_t^{x,i}\left[\int_{t+\epsilon}^Te^{-\rho_i(s-t)}
[U(F_2(s,\bar{X}(s),J(s)))-U(c(s))]\,ds\right]\nonumber\\
&+&\mathbb{E}_t^{x,i}\left[e^{-\rho_i(T-t)}(U(\bar{X}(T))-U({X}(T)))\right].
\end{eqnarray}
In the light of inequalities \eqref{089} and \eqref{189} and
Dominated Convergence Theorem
\begin{eqnarray}
{\lim_{\epsilon\downarrow 0}}
\frac{\mathbb{E}_t^{x,i}\left[\int_{t}^{t+\epsilon}e^{-\rho_i(s-t)}
[U(F_2(s,\bar{X}(s),J(s)))-U(c(s))]\,ds\right]}{\epsilon}=
U(F_2(t,x,i)-U(c(t)).\nonumber \end{eqnarray}
In the light of \eqref{yy} it follows that
\begin{eqnarray}\label{7.19}
&&\mathbb{E}_t^{x,i}\left[\int_{t+\epsilon}^Te^{-\rho_i(s-t)}
[U(F_2(s,\bar{X}(s),J(s)))-U(c(s))]\,ds\right]+\mathbb{E}_t^{x,i}
\left[e^{-\rho_i(T-t)}(U(\bar{X}(T))-U({X}(T)))\right]\nonumber\\
&=&\mathbb{E}_t^{x,i}\left
[v(t+\epsilon,\bar{X}(t+\epsilon),J(t+\epsilon))
-v(t+\epsilon,{X}(t+\epsilon),J(t+\epsilon))\right]\nonumber\\
&+&\mathbb{E}_t^{x,i}\left[\mathbb{E}[\int_{t+\epsilon}^T
(e^{-\rho_i(s-t)}-e^{-\rho_{J(t+\varepsilon)}(s-t)})
[U(F_2(s,\bar{X}(s),J(s)))-U(F_2(s,{X}(s),J(s)))]\,ds\right]\nonumber\\
&+&\mathbb{E}_t^{x,i}\left[\mathbb{E}[
(e^{-\rho_i(T-t)}-e^{-\rho_{J(t+\varepsilon)}(T-t+\epsilon)})
[U(\bar{X}(T))-U({X}(T))]\right]
\end{eqnarray}

 In the light of inequalities \eqref{089}, \eqref{189} and
Dominated Convergence Theorem it follows that

$${\lim_{\epsilon\downarrow 0}}\frac{\mathbb{E}_t^{x,i}\left[
(e^{-\rho_i(T-t)}-e^{-\rho_{J(t+\varepsilon)}(T-t+\epsilon)})
[U(\bar{X}(T))-U({X}(T))]\right]}{\varepsilon}$$

$$=\lambda_{ij}
(e^{-\rho_i(T-t)}-e^{-\rho_{j}(T-t)})\mathbb{E}_t^{x,i}\left[
[U(\bar{X}(T))-U({X}(T))]\right]=0.$$
By the same token one can get that
$${\lim_{\epsilon\downarrow 0}}\frac{
\mathbb{E}_t^{x,i}\left[\int_{t+\epsilon}^T
(e^{-\rho_i(s-t)}-e^{-\rho_{J(t+\varepsilon)}(s-t)})
[U(F_2(s,\bar{X}(s),J(s)))-U(F_2(s,{X}(s),J(s)))]\,ds\right]}{\epsilon}=0.$$
It\^{o}'s formula yields
\begin{eqnarray*}
&&\mathbb{E}_t^{x,i}\left[v(t+\epsilon,\bar{X}(t+\epsilon),J(t+\epsilon))
-v(t+\epsilon,{X}(t+\epsilon),J(t+\epsilon))\right]\nonumber\\
&=&\mathbb{E}_t^{x,i}\int_t^{t+\epsilon}[\overline{\Gamma}
v(s,\bar{X}(s),J(s))- U(F_{2}(s,\bar{X}(s),J(s))]ds -\mathbb{E}_t^{x,i}\int_t^{t+\epsilon}[\Gamma
v(s,{X}(s),J(s))-U(c(s))]ds,
\end{eqnarray*}

Therefore

\begin{eqnarray}
{\lim_{\epsilon\downarrow 0}}\frac{J(t,x,i,F_{1},F_{2})-J(t,x,i,
\pi_{\epsilon},c_{\epsilon})}{\epsilon}
= [\overline{\Gamma} v(t,x,i)-{\Gamma} v(t,x,i)]\geq0,\nonumber.
\end{eqnarray}
the inequality follows from \eqref{21}.

\begin{flushright}
$\square$
\end{flushright}

\subsection{C Proof of Lemma \ref{4.4}}
 We know that $h(t,i)$ solve the ODE system
$$\frac{\partial h}{\partial t}(t,{\bf{0}})+(\lambda_{11}-\rho_{\bf{0}}) h(t,{\bf{0}})+
\lambda_{12}h(t,{\bf{1}})+1=0,$$
$$\frac{\partial h}{\partial t}(t,{\bf{1}})+(\lambda_{22}-\rho_{\bf{1}}) h(t,{\bf{1}})+
\lambda_{21}h(t,{\bf{0}})+1=0.$$ Let us define $h(t)\triangleq h(t,{\bf{0}})-h(t,{\bf{1}})$, then
 $h$ satisfies the following ODE:
$$h'(t)+(\lambda_{11}+\lambda_{22}-\rho_{\bf{0}})h(t)-(\rho_{\bf{0}}-\rho_{\bf{1}})h(t,{\bf{1}})=0$$
Under the assumption $\rho_{\bf{0}}>\rho_{\bf{1}},$ one gets
$$h'(t)+(\lambda_{11}+\lambda_{22}-\rho_{\bf{0}})h(t)>0.$$
From here it follows that $h(t)<0,$ on $[0,T)$ (since $h(T)=0$). Hence, $h(t,{\bf{0}})<h(t,{\bf{1}}).$ Consequently $F_2(t,x,{\bf{0}})>F_2(t,x,{\bf{1}}).$

\begin{flushright}
$\square$
\end{flushright}

\subsection{ D Proof of Lemma \ref{merton}}
Here we are in the case of $\rho_{{\bf{0}}}=\rho_{{\bf{1}}}=\rho.$ Further, for simplicity we assume
constant market coefficients $\alpha$ and $\sigma.$ The optimal consumption rate is
\begin{equation}\label{c*}
\bar{C}(t)=\frac{\bar{c}(t)}{X(t)}=\frac{\eta}{(1+(\eta-1)e^{\eta(t-T)})},
\end{equation}
where $\eta\triangleq\frac{\rho-\gamma[\frac{(\alpha-r)^2}{2\sigma^2(1-\gamma)}+r]}{1-\gamma}.$ It follows from direct computations that
\begin{eqnarray}\label{pc}
\frac{\partial \bar{C}}{\partial
\rho}=-\frac{1}{\bar{C}^{2}(t)}[-\frac{1}{(1-\gamma)\eta^2}+e^{\eta(t-T)}\frac{1+\eta(\eta-1)(t-T)}{\eta^2(1-\gamma)}]
\end{eqnarray}
If $1+\eta(\eta-1)(t-T)\leq0$ the claim follows. Otherwise, since
$$e^{\eta(T-t)}>1+\eta(T-t)>1+\eta(T-t)+\eta^2(t-T),$$
it follows that
\begin{equation}\label{e}
e^{\eta(t-T)}<\frac{1}{1+\eta(\eta-1)(t-T)},
\end{equation}
whence the result.

\begin{flushright}
$\square$
\end{flushright}


\begin{thebibliography}{99}




\bibitem{Das}
{\sc Das M.} (2003) Optimal growth with decreasing marginal
impatience, {\em Journal of Economic Dynamics \& Control}, {\bf 27}, 1881-1898.

\bibitem{EkeLaz}
{\sc Ekeland, I. and Lazrak, A.} (2010), The golden rule when preferences are time inconsistent, {\em Mathematics and Financial Economics}, {\bf{4}}, 1-27.

\bibitem{EkePir}
{\sc Ekeland, I. and Pirvu, T. A.} (2008), Investment and consumption without commitment, {\em Mathematics and Financial Economics}, {\bf{2}}, 57-86.

\bibitem{EkePirMbo}
{\sc Ekeland, I., Mbodji, O., and Pirvu, T. A.} (2010)
Time consistent portfolio
management,  {\em http://arxiv.org/pdf/1008.3407 }.

\bibitem{Eli} \textsc{Elliott, R. J., Chan, L., and Siu, K., T.} (2005)
{Option pricing and Esscher transform
under regime switching}, {\em Annals of Finance} {\bf 1}, 423-432.

\bibitem{Eps}
{\sc Epstein, L. G. and Hynes, J. A.} (1983) The rate of time preference and dynamic economic analysis, {\em Journal of Political Economy}, {\bf 91}, 611-635

\bibitem{Guo} \textsc{Guo, X.} (2001) {Information and option pricing} , {\em Quantitative Finance} {\bf 1}, 38-44.


\bibitem{Koop}
{\sc Koopmans, T.C.} (1960) Stationary utility and impatience, {\em Econometrica} {\bf
28}, 287-309.

\bibitem{KruSm}
{\sc Krusell, P and Smith, A.} (2003)  Consumption and Savings
Decisions with Quasi-geometric Discounting, {\em Econometrica} {\bf
71}, 365-375.


\bibitem{Lai}
{\sc Laibson, D.} (1997)  Golden Eggs and Hyperbolic Discounting
{\em The Quarterly Journal of Economics}, {\bf 112}, 443-477.

\bibitem{Mer69}
{\sc Merton, R. C.} (1969) Lifetime portfolio selection under
uncertainty: the continuous-time case, {\em Rev. Econom. Statist.}
{\bf 51}, 247-257.

\bibitem{Mer71}
{\sc  Merton, R.C.} ( 1971)
\newblock Optimum consumption and portfolio rules in a conitinuous-time model.
\newblock {\em J. Economic Theory}, {\bf 3}, 373-413.

\bibitem{STRO}
{\sc Strotz R.} (1955)  Myopia and inconsistency in Dynamic Utility
Maximization, {\em Rev. Financial Stud.}, {\bf 23}, 165-180.


\bibitem{Cad}
{\sc Sotomayor, L. R. and Cadenillas, A.} (2009) Explicit Solutions of Consumption-Investment Problems in Financial Markets with Regime Switching, {\em Mathematical Finance}, {\bf 9} (2), 251-279.



\bibitem{Tak}
{\sc Takashi K.} (2000) Increasing Marginal Impatience
and Intertemporal Substitution {\em Journal of Economics}, {\bf 72}, 67-79.

\bibitem{Tosh}
{\sc Toshihiko M.} (2009) A note on cyclical discount factors and labor market volatility. {\em Preprint}.

\bibitem{Uza}
{\sc Uzawa, H.} (1968) Time preference, the consumption function, and optimum holdings {\em
In: Wolfe, J.N. (Ed.), Value, Capital and Growth: Papers in Honour of Sir John Hicks. Edinburgh University Press, Edinburgh.}







\end{thebibliography}
\end{document}